\documentclass[11pt]{article}

\usepackage{amsmath}
\usepackage{amsthm}
\usepackage{amssymb}
\usepackage{amscd}
\usepackage[all]{xy}

\title{On the canonical threefolds with strictly nef anticanonical divisors}
\author{Hokuto Uehara}
\date{}

\theoremstyle{plain}
\newtheorem{prop}{Proposition}[section]
\newtheorem{lem}[prop]{Lemma}
\newtheorem{thm}[prop]{Theorem}

\theoremstyle{definition}
\newtheorem{defn}[prop]{Definition}

\theoremstyle{remark}
\newtheorem{rem}[prop]{Remark}

\newcommand{\Alb}{\mathrm{Alb}}
\newcommand{\alb}{\mathrm{alb}}

\begin{document}

\maketitle

\begin{abstract}
The purpose of this article is to give a criterion for canonical 3-folds to be a $\mathbb{Q}$-Fano 3-folds.
\end{abstract}

\setcounter{section}{-1}
\section{Introduction}

\indent

We are working in the category of projective varieties defined over the field of complex numbers.
 
 We study a numerical criterion for canonical 3-folds to be $\mathbb{Q}$-Fano 3-folds. A $\mathbb{Q}$-Fano 3-fold is defined in this paper to be a normal projective variety of dimension 3 having only canonical singularities with its anticanonical divisor being ample.    
 
  First of all, we state the following famous criterion for ampleness of divisors, so-called Kleiman's criterion.

\begin{thm}\textbf{\textup{\cite{Kleiman}}}\label{Kleiman}
Let $X$ be a projective variety and $D$ a Cartier divisor. Then $D$ is ample if and only if $D\cdot z>0$ for all $z\in \overline{NE}(X)\backslash \{0\}$.
\end{thm}

Next we introduce the notion, \textit{strictly nef}.   
\begin{defn}
 Let $X$ be a normal projective variety and $D$ a $\mathbb{Q}$-Cartier Weil
divisor. $D$ is \textit{strictly nef}, if $D\cdot C> 0$ for all irreducible
curves $C$.
\end{defn}
 We are interested in how far strictly nef divisors are from being ample. Mumford gave an example of a strictly nef non-ample divisor on a surface (see \cite{Ha70}, p56). In words in Theorem \ref{Kleiman}, $D\cdot z>0$ for all $z\in NE(X)\backslash \{0\}$ does not imply that $D$ is ample. But in the case $D=-K_{X}$, we have the following results by H. Maeda, F. Serrano and K. Matsuki. 

\begin{thm} \label{Serrano3}
Let $X$ be a $\mathbb{Q}$-Gorenstein normal projective variety with strictly nef $-K_X$. Assume that one of the following conditions holds. 
\renewcommand{\labelenumi}{(\roman{enumi})}
\begin{enumerate}
\item \textup{\cite{Maeda}} $X$ is smooth and $\dim X=2$.
\item \textup{\cite{Serrano}} $X$ is smooth and $\dim X=3$.
\item \textup{\cite{Matsuki}} $X$ is a canonical 3-fold and $\kappa(-K_X)\ge0$.
\end{enumerate}
Then $-K_X$ is ample. 
\end{thm}
 Furthermore in more general aspect, we have
\begin{thm} \textbf{\textup{\cite{Serrano}}}\label{Serrano2}
Let $X$ be an irreducible reduced projective Gorenstein surface and L a strictly nef Cartier divisor. Then $mL+K_X$ is ample for a sufficiently large integer $m$. 
\end{thm}

 Our purpose in this paper is to prove the following theorem. We will prove it with the aid of various results including Theorem \ref{Serrano3} and \ref{Serrano2}. \\

\noindent
\textbf{Main Theorem (=Theorem \ref{q=0}).} {\itshape{Let $X$ be a canonical 3-fold with strictly nef $-K_X$. Then $-K_{X}$ is ample.}} \\

 The proof rests mainly on sophisticated techniques of three dimensional birational geometry developed in the last two decades, essentially by the Japanese School. We now give an outline of this article. Let $X$ be as in the Main Theorem. Let $\alpha :X \to \Alb(X)$ be the Albanese morphism and $\alb (X)$ denotes $\dim\alpha (X)$. We will show in \S 1 that the case $\alb(X)=1$ never occurs and in \S 2 that the case $\alb(X)=2$ never occurs by the extremal ray theory. Since $\kappa(X)=-\infty$, we can check $\alb(X)\ne 3$ easily (Remark \ref{q==0}). After we know them, it is enough to show the ampleness of $-K_{X}$ with the assumption $q(X)=0$ in order to prove the Main Theorem. It will be done in \S 3. \\

\noindent
\textit{Acknowledgment.} 
I express my hearty thanks to Professor Yujiro Kawamata, my thesis advisor, for giving me useful comments and invaluable advice. I also thank Dr. Hiromichi Takagi and Dr. Tatsuhiro Minagawa for stimulating discussions on this paper. Finally I am grateful to Professor Hajime Kaji for giving me the opportunity to talk about this subject in his master seminar while preparing this paper. \\

\noindent
\textit{Notation and Convention.} 

\noindent
(i) We work over the field of complex numbers.

\noindent
(ii) A terminal (resp. canonical) 3-fold means a normal projective 3-dimensional variety with at most terminal (resp. canonical) singularities.

\noindent
(iii) A surjective morphism between normal projective varieties is an algebraic fiber space if it has connected fibers.

\noindent
(iv) An algebraic fiber space $\varphi:X\to Y$ between normal projective varieties is an extremal contraction if $-K_{X}$ is $\varphi$-ample and if the Picard numbers satisfy $\rho(X)=\rho(Y)+1$.

\noindent
(v) A birational morphism between normal projective varieties is small if it is isomophic in codimension 1.
 
\noindent
(vi) We always denote numerical equivalence by $\equiv$.
 
\noindent
(vii) Let $X$ be as in the Main Theorem. We always denote the Albanese morphism of $X$ by $\alpha$ in this paper. Define $\alb (X):= \dim\alpha (X)$.\\

In the whole paper we will freely use the results from classification theory and Mori theory. Refer to \cite{KMM} and \cite{KM}. 
  
\section{Proof of $\boldsymbol{\alb(X)\ne1}$}

\indent
 
 A normal variety $X$ is said to have only \textit{terminal singularities} (resp. \textit{canonical singularities}) if the following conditions are satisfied: 

(i) The canonical divisor $K_X$ is $\mathbb{Q}$-Cartier, i.e. $mK_X$ is a Cartier divisor for some positive integer $m$.

(ii) For a resolution of singularities $f:Y\to X$ we have
\[K_Y=f^* K_X +\sum{a_i E_i}\]
with $a_i>0$ (resp. $a_i\ge0$) for all $i$, where the $E_i$'s run through all the prime divisors on $Y$ which are exceptional with respect to $f$. In paticular if $X$ is a surface, $X$ has only terminal singularities (resp. canonical singularities) is equivalent to that $X$ is nonsingular (resp. has only Du Val singularities, i.e. rational double points). 

 Let $X$ be a normal projective variety with at most canonical singularities. Then $X$ has rational singularities (\cite{Elkik}) so the Albanese map of $X$ is a morphism (cf. Lemma 8.1 in \cite{K85}). 
 
 Assume in addition that $\dim X=3$. Because a general hyperplane section of $X$ has at most Du Val singularities, the locus of non-cDV (compound Du Val) points is a finite set. See also Corollary 5.40 in \cite{KM}.

 We enjoy the benefit of the following result repeatedly in this paper.
 
\begin{prop} \textbf{\textup{\cite{KeMM}}}\label{KeMM}
Let $X$ be a $\mathbb{Q}$-factorial terminal 3-fold such that $-K_{X}$ is nef, but not big. Then $c_{1}(X)c_{2}(X) \ge 0$. 
\end{prop}
 
\begin{rem} \label{KeMM rem}
Taking a $\mathbb{Q}$-factorialization allows us to remove the assumption $\mathbb{Q}$-factoriality in Proposition \ref{KeMM}. For the definition of $\mathbb{Q}$-factorialization, see Remark \ref{Q-fac}. By \cite{K86}, there exists a nonnegative rational number $A$ such that $\chi(\mathcal{O}_X)=\frac{c_{1}c_{2}}{24}+A$, where $A=0$ if and only if $X$ is Gorenstein. Thus assume in addition that $\chi(\mathcal{O}_X)=0$ in the proposition above. Then $X$ is Gorenstein.
\end{rem} 
 
\begin{lem}\textbf{\textup{\cite{Serrano}}}\label{serrano-1}
Let $Y$ be a smooth projective variety, $C$ a smooth projective curve of positive genus, and $g:Y\to C$ a morphism with connected fibers. We denote by $F_t$ the fiber over $t\in C$. Let $D$ be a divisor on $Y$ and pick a paticular fiber $F_0$, $0\in C$. If $H^0(F_0,D|_{F_0})\ne 0$, and $H^1(Y,D+F_t-F_0)=0$ for all $t\in C$, then $D$ is algebraically equivalent to an effective divisor.  
\end{lem}  

 Using the idea of \cite{Serrano}, we can show the following proposition rather easily. The key fact of the proof is that $\mathbb{Q}$-factorial Gorenstein terminal 3-folds are factorial (Lemma 5.1 in \cite{K88}).       
     
\begin{prop}\label{q=1}
  Let $X$ be a $\mathbb{Q}$-factorial terminal 3-fold with strictly nef $-K_{X}$. Then $q(X)=1$ never occurs.
\end{prop}

\begin{proof}
Assume that $q(X)=1$. Then $\alpha$ is an algebraic fiber space and $C:=\Alb(X)$ is an elliptic curve. Let $F$ be a general fiber of $\alpha$. Then $F$ is a smooth projective surface, so $-K_F=-K_X|_{F}$ is ample (Theorem \ref{Serrano3}). Hence $R^i \alpha _* \mathcal{O}_X=0$ for $i>0$ by Kawamata-Viehweg Vanishing Theorem (cf. \cite{KMM}, Theorem 1.2.5) and so $h^{i}(\mathcal{O}_X)=h^i(\mathcal{O}_C)$ for all $i$. Then applying Remark \ref {KeMM rem}, we know that $X$ is a $\mathbb{Q}$-factorial Gorenstein terminal 3-fold.  

 Let $f: Y\to X$ be a resolution and put $g:=\alpha \circ f$.\\
 
\underline{Claim.} $-mK_X \equiv E$ for some effective divisor $E$ and some $m>0$.\\

\noindent \textit{Proof of the claim}. Set $X_t$ (resp. $Y_t$) the fiber of $\alpha$ (resp. $g$) over $t\in C$. For a sufficiently large integer $m$, 
\begin{equation}
h^{0}(F,-mf^{*}K_{X}|_F)=h^{0}(F,-mK_F)\ne 0
\end{equation} 
 by the ampleness of $-K_F$. Fix this $m$. We may assume $h^{0}(-mf^{*}K_X+F-Y_t)=h^{0}(-mK_X+F-X_t)=0$. Let $S$ be a sufficiently ample line bundle on $X$. We consider the following exact sequence.
\begin{align*}
\cdots \to H^{1}((m+1)K_{X}-F+X_{t}-S) &\to H^{1}((m+1)K_{X}-F+X_t) \\
&\to H^{1}((m+1)K_{X}-F+X_{t}|_S) \to\cdots
\end{align*}
Since $(-K_{X})^{2}\not \equiv 0$, $-(m+1)K_X+F-X_{t}|_S$ is nef and big. Thus 
\[h^{1}((m+1)K_{X}-F+X_{t}|_S)=h^{1}(-(m+1)K_{X}+F-X_{t}|_S+K_S)=0.\] 
 Therefore looking at the exact sequence above,
\[h^{2}(-mf^{*}K_X+F-Y_t)=h^{2}(-mK_X+F-X_t)=h^{1}((m+1)K_X-F+X_t)=0.\]
Thus
\begin{equation}
\quad -h^{1}(-mf^{*}K_X+F-Y_t)=\chi(-mK_X+F-X_t)=(2m+1)\chi(\mathcal{O}_X)=0 
\end{equation}  
by Riemann-Roch Theorem (cf. \cite{Reid87}). Then we can apply Lemma \ref{serrano-1} because of (1) and (2). Thus we know that $-mf^{*}K_X$ is algebraically equivalent to an effective divisor. The algebraically equivalence of cycles is preserved by proper push-forward (see \cite{Fu97}, Proposition 10.3). Then we obtain the claim. \\

 We can write $-K_{X}\equiv\displaystyle\sum a_{i}Z_i$ for some $a_{i}\in \mathbb{Q}_{>0}$ and some prime divisors $Z_i$ by the claim. Since $-K_{X}$ is nef but not big by $q(X)\ne 0$, $(-K_X)^3=0$ and so ${(-K_X)}^{2}Z_i=0$ for all $i$. Because $Z_i$ is Cartier by Lemma 5.1 in \cite{K88}, $Z_i$ is a Gorenstein surface. Hence $K_{Z_i}-nK_{X}|_{Z_i}$ is ample for $n\gg 0$ (Theorem \ref{Serrano2}), then 
 \[0<-K_{X}|_{Z_i}(K_{Z_i}-nK_{X}|_{Z_i})=-K_{X}Z_{i}(K_{X}+Z_{i}-nK_{X})=-K_{X}{Z_i}^2 .\]
 The absurdity comes from 
 \[ 0=(-K_{X})^2Z_{1}=(-K_X)(\sum {a_i Z_i})Z_{1}=a_{1}(-K_X ){Z_1}^2+(-K_X)\sum_{i\ne1}a_{i}Z_{1}Z_{i}>0.\] 
\end{proof}

\begin{lem}\label{q>1}
Let $X$ be a canonical 3-fold with strictly nef $-K_{X}$. Assume in addition that $\alb(X)=1$. Then $q(X)=1$.
\end{lem}

\begin{proof}
 Note that $\alpha$ is an algebraic fiber space in the case $\alb(X)=1$ when we regard $\alpha$ as a morphism from $X$ to $\alpha(X)$. Assume $q(X)>1$. Denote a terminalization of $X$ by $f:Y\to X$ (see Remark \ref{term1}). Let $F$ be a general fiber of $\alpha$. Then $F':=f^*(F)$ is smooth and connected, so $F$ is an irreducible reduced Gorenstein surface. Therefore $-K_{F}=-K_X|_{F}$ is ample by Theorem \ref{Serrano2} and we obtain $R^i\alpha _*\mathcal{O}_X=0$ for all $i>0$ by relative Kawamata-Viehweg Vanishing Theorem. Then $h^{i}(\mathcal{O}_Y)=h^i(\mathcal{O}_X)=h^{i}(\mathcal{O}_{\alpha(X)})$ for all $i$, and so $\chi(\mathcal{O}_Y)=\chi(\mathcal{O}_{\alpha(X)})<0$. This lies in contradiction with Proposition \ref{KeMM}.
 \end{proof}

 The following definition is in \cite{PS}. It is useful, for example, when we investigate a $\mathbb{Q}$-factorial terminal 3-fold with nef anticanonical divisor.
\begin{defn}
Let $X$ be a normal projective variety and $D$ a $\mathbb{Q}$-Cartier Weil
divisor on $X$.
Then $D$ is called \textit{almost nef}, if there are at most finitely many
rational curves $C_i$, $1\leq i\leq r$, such that $D\cdot C\geq 0$
for all curves $C\not= C_i$.
\end{defn}

\begin{prop}\textbf{\textup{\cite{PS}}}\label{prop:PS1}
Let $X$ be a $\mathbb{Q}$-factorial terminal 3-fold with almost nef
$-K_{X}$.
Assume that one of the following conditions holds.
\renewcommand{\labelenumi}{(\roman{enumi})}
\begin{enumerate}
\item There is an extremal contraction $\varphi :X\to C$ to the elliptic
curve $C$. 
\item $q(X)=1$ and there exists an extremal contraction $\varphi :X\to W$ 
to the surface $W$ .
\item $\alb(X)=2$ and there exists an extremal contraction $\varphi :X\to W$
to the surface $W$ .
\end{enumerate}
Then the Albanese morphism $\alpha:X\to \Alb(X)$ is smooth and $-K_{X}$ is nef.
\end{prop}

 We will use Proposition \ref{prop:PS1} under the condition (iii) only.  
 
\begin{prop}\textbf{\textup{\cite{PS}}}\label{prop:PS}
Let $X$ be a $\mathbb{Q}$-factorial terminal 3-fold with almost nef $-K_{X}$.
\renewcommand{\labelenumi}{(\roman{enumi})}
\begin{enumerate}
\item Let $\varphi :X\to W$ be a divisorial contraction or a flip. Then $-K_{W}$ is almost nef. 
\item Let $\varphi :X\to W$ be an extremal contraction to the surface $W$. Then $-(4K_W+\Delta )$ is almost nef \textup{(}see the definition of $\Delta$ in \textup{\cite{PS}}, which is an effective divisor on $W$\textup{)}.
\end{enumerate}
\end{prop}

\begin{rem}
Let $X$ be a $\mathbb{Q}$-factorial terminal 3-fold with nef $-K_{X}$ and $\varphi :X\to W$ a divisorial contraction. This does not imply that $-K_W$ is \textit{nef} (see Proposition 3.3 in \cite{DPS93}). As in \cite{PS}, this is the reason why they introduce the notion \textit{almost nef}.  
\end{rem}

\begin{rem}\label{term1}
Let $X$ be a 3-fold with canonical singularities. According to \cite{Reid83}, there exists a partial resolution $f:Y\to X$ such that (i) $Y$ has terminal singularities, and (ii) $f$ is crepant and projective. We call $f$ and $Y$ terminalization of $X$.
\end{rem}

We have more information on terminalization as follows. We denote the number of crepant divisors over a canonical 3-fold $X$ by $e(X)$.

\begin{prop}\textbf{\textup{\cite{K88}}}, \textbf{\textup{\cite{Kollar89}}}\label{term2}
Let $X$ be a normal projective 3-dimensional variety having at most $\mathbb{Q}$-factorial canonical singularities. Assume that $e(X)>0$. Then there exist a normal projective variety $X_1$ having at most $\mathbb{Q}$-factorial canonical singularities and a projective birational morphism $g:X_1\to X$ such that the following conditions are satisfied:
\renewcommand{\labelenumi}{(\roman{enumi})}
\begin{enumerate}
\item $g$ is crepant; 
\item The exceptional locus of $g$ is a prime divisor;
\item $\rho(X_{1}/X)=1$ and $e(X_1)=e(X)-1$.
\end{enumerate}  
\end{prop} 

Using Proposition \ref{term2} iterately, for a given $\mathbb{Q}$-factorial canonical 3-fold $X$, we obtain its terminalization because $e(X)<\infty$. 

\begin{rem}\label{Q-fac}
Recall that a $\mathbb{Q}$-factorialization of a canonical 3-fold $X$ means a birational projective morphism $f:Y\to X$ such that (i) $Y$ has at most $\mathbb{Q}$-factorial canonical singularities, and (ii) $f$ is small. Yujiro Kawamata proved in \cite{K88} the existence of $\mathbb{Q}$-factorialization of canonical 3-folds. Note that if $f(C)$ is a point for an irreducible curve $C$ on $Y$, then $C=\mathbb{P}^1$ because $R^i f_* \mathcal{O}_Y =0$ for $i>0$.
\end{rem}
   
The following proposition is the main result of this section.

\begin{prop}\label{alb=1}
Let $X$ be a canonical 3-fold with strictly nef $-K_{X}$. Then $\alb(X)\ne1$. 
\end{prop}

\begin{proof}
 Assume that $\alb(X)=1$. Then we know that $q(X)=1$ because of Lemma \ref{q>1} and $\Alb(X)$ is an elliptic curve. Let $g_0:Y_0\to X$ be a $\mathbb{Q}$-factorialization. Then there exists a sequence
\[ Y:=Y_n\to Y_{n-1}\to\dots\to Y_1 \to Y_0\to X. \]
such that for $i\ge1$ each morphism $g_i:Y_i\to Y_{i-1}$ is constructed by Proposition \ref{term2} and $Y:=Y_n$ has at most $\mathbb{Q}$-factorial terminal singulalities (if $e(Y_0)=0$, $n=0$). Put $f:=g_0\circ g_1\circ \dots \circ g_n$. Since $K_{Y}$ is not nef, there exists an extremal contraction $\varphi :Y\to W$. Since $h^{i}(\mathcal{O}_Y)=h^{i}(\mathcal{O}_{\Alb(X)})$ for $i\ge 0$ as in the proof of Lemma \ref{q>1}, we obtain that $Y$ is Gorenstein (Remark \ref{KeMM rem}). Furthermore we know that $\dim W\ge1$ easily. Let $\beta:W\to \Alb(W)(=\Alb(X))$ be the Albanese morphism.\\
 
\textbf{\underline{Case: $\dim W=1$}} 
 
 In this case, $W$ is an elliptic curve and so $\beta$ is an isomorphism. Proposition \ref{q=1} says that $X$ is not $\mathbb{Q}$-factorial terminal, so there exists an irreducible curve $C$ such that $f(C)$ is a point. Then $K_Y\cdot C=0$ and so $\varphi (C)$ is a curve. This lies in contradiction with $\alpha\circ f=\varphi$. \\
 
\textbf{\underline{Case: $\dim W=2$}}    
 
 When $\kappa(W)\ge 0$, $K_{W}\cdot H\ge 0$ for some ample divisor $H$. But because $-(4K_{W}+\Delta)$ is almost nef (Proposition \ref{prop:PS}), $-K_{W}\cdot H\ge0$. Therefore we know $K_{W}\equiv0$. Combining this with $h^{i}(\mathcal{O}_W)=h^{i}(\mathcal{O}_Y)=h^{i}(\mathcal{O}_{\Alb(X)})$ for all $i$, we know that $W$ is a hyperelliptic surface. If $X$ has only terminal singularities, $f$ is just a $\mathbb{Q}$-factorialization so it is small. $X$ is not $\mathbb{Q}$-factorial terminal by Proposition \ref{q=1}, then there exists a rational curve $C$ on $Y$ such that $K_Y\cdot C=0$, and so $\varphi(C)$ is a rational curve on a hyperelliptic surface $W$. This is absurd. Thus $e(X)>0$, that is to say, $X$ is not terminal and there exists a morphism $g_n:Y\to Y_{n-1}$ such that $g_n$ is constructed by Proposition \ref{term2} and $f$ factors through $g_n$. Put $D$ the exceptional divisor of $g_n$. If $\dim g_n(D)=0$, then $D$ is uniruled by the Subadjunction Lemma (Lemma 5.1.9 !
of \cite{KMM}). But this lies in contradiction with $\varphi (D)=W$. If $\dim g_n(D)=1$, then there exists an irreducible curve $C$ on $Y$ such that $g_n(C)$ is a point. Since $\dim {g_n}^{-1}(g_n(C))=1$, $C$ is a smooth rational curve. But this derives contradiction as above.   
  
 Assume that $\kappa(W)=-\infty$ below. We obtain that $-K_Y \equiv \displaystyle\sum a_{i}Z_i$ for some $a_{i}\in \mathbb{Q}_{>0}$ and some prime divisors $Z_i$ as in the proof of Proposition \ref{q=1}. Since $(-K_{Y})^3=0$, $(-K_{Y})^2\cdot Z_{i}=0$ and there exists an irreducible curve $C_{i}'\subset Z_{i}$ such that $f(C_{i}')$ is a point for all $i$. In fact, if not, $-K_{Y}|_{Z_{i_0}}$ is strictly nef for some $i_0$. Moreover we know by Theorem \ref{Serrano2} that 
 \[0<-K_{Y}|_{Z_{i_0}}(K_{Z_{i_0}}-nK_{Y}|_{Z_{i_0}})=-K_{Y}Z_{i_0}(K_{Y}+Z_{i_0}-nK_{Y})=-K_{Y}Z_{i_0}^2 .\]
 for $n\gg0$. Then we can derive a contradiction as follows:
 \[ 0=(-K_{Y})^2 Z_{i_0}=(-K_Y)(\sum {a_i Z_i})Z_{i_0}=a_{i_0}(-K_Y)(Z_{i_0})^2+(-K_Y)\sum_{i\ne i_0}a_{i}Z_{i_0}Z_{i}>0.\] 
   
 Let $\gamma :Y\to \Alb(Y)=(\Alb(X))$ be the Albanese morphism. Note that $\alpha\circ f=\gamma=\beta\circ \varphi$.
\[
{\xymatrix{
X \ar[d]_{\alpha} & Y \ar[l]_{f} \ar[dl]_{\gamma} \ar[d]^{\varphi} \\
\Alb(X) & W \ar[l]^{\beta} }}
\]

 \underline{Claim.} $\varphi(Z_{i_0})=W$ for some $i_{0}$.\\
 
  \noindent \textit{Proof of the claim}. If $\gamma(Z_{i})$ is a point for all $i$, then there exists an irreducible curve $C$ on $Y$ such that $f(C)$ is 1-dimensional and $C$ does not meet any $Z_{i}$. This implies
 \[0<-K_{Y}\cdot C=\sum a_{i}Z_{i}\cdot C=0\]
 which is contradiction. Therefore we know that $\gamma (Z_{i_0})=\Alb(X)$ for some $i_{0}$, then $\varphi(Z_{i_0})=W$. In fact, assume to the contrary that $\varphi(Z_{i_0})$ is a curve. Now that $\varphi(C_{i_0}')$ is a curve, $\Alb(X)=\gamma (Z_{i_0})=\beta\circ\varphi(Z_{i_0})=\beta\circ\varphi(C_{i_0}')$. But this is absurd because $\beta\circ\varphi(C_{i_0}')$ is a point. \\  
 
 Using the claim above, we prove the following.\\ 
 
 \underline{Claim.} Assume that $-(4K_{W}+\Delta)\cdot \tilde{C}=0$ for an irreducible curve $\tilde{C}$ on $W$. Then $\tilde{C}$ is a fiber of $\beta$. \\

\noindent \textit{Proof of the claim}. If $f(C)$ is a curve for all curves $C\subset \varphi ^{-1}(\tilde{C})$, then $Z_{i}\not \subset \varphi ^{-1}(\tilde{C})$ for all $i$ by $C_{i}'\subset Z_{i}$. From the well known formula $(K_Y)^2\cdot \varphi ^{-1}(\tilde{C})=-(4K_{W}+\Delta)\cdot \tilde{C}$ (see the proof of Lemma 1.6 in \cite{PS}), we get
 \[-K_{Y}\cdot (\sum a_{i}Z_{i})\cdot \varphi ^{-1}(\tilde{C}) =0.\] 
Then $\dim Z_{i}\cap \varphi ^{-1}(\tilde{C})\le 0$ for all $i$. This lies in contradiction with $\varphi(Z_{i_0})=W$. Then we know that there exists a curve $C\subset \varphi ^{-1}(\tilde{C})$ such that $f(C)$ is a point, and so $\tilde{C}=\varphi (C)$ is a fiber of $\beta$. \\

 It follows from Proposition 1.7 in \cite{PS} that $W$ is a $\mathbb{P}^1$-bundle over $\Alb(X)$ with nef $-K_{W}$ and almost nef $-(4K_{W}+\Delta)$. Then we can write $W=\mathbb{P}^{1}(V)$ for a normalized sheaf $V$ on $C$, i.e. $H^{0}(V)\neq0$, but $H^{0}(V\otimes \mathcal{L})=0$ for all line bundles $\mathcal{L}$ with $\deg\mathcal{L}<0$ (see [7] about the treatment and the terminology of geometrically ruled surfaces around here ). Set $e:=-c_1(V)$. We separate into these two cases: $e=0$ or $-1$. Let $F$ be a fiber of $\beta$ and $C_{0}$ a canonical section satisfying $C_{0}^2=-e$.
 
 In the case $e=-1$, there exists an irreducible curve $C_{1}$ on $W$ such that $C_{1}$ is numerically proportional to $2C_{0}-F$. Since $C_{1}$ is nef,
\[0\le-(4K_{W}+\Delta)\cdot C_{1}=-\Delta\cdot C_{1}\le 0.\]
Thus 
\[ -(4K_{W}+\Delta)\cdot C_{1}=0.\]
This contradicts the claim above.

In the case $e=0$, $C_{0}$ is nef. Then we obtain 
\[ 0\le-(4K_{W}+\Delta)\cdot C_{0}=-\Delta\cdot C_{0}\le 0.\]
This also yields contradiction. \\

 \textbf{\underline{Case: $\dim W=3$}} 
   
 Since $Y$ is Gorenstein, $\varphi$ is a divisorial contraction (cf. \cite{Benveniste}). Let $E$ be the exceptional divisor. If $\varphi(E)$ is a point, then we have the following equation 
\[ K_Y=\varphi ^{*}K_W+aE \] 
for some $a\in\mathbb{Q}_{>0}$. Thus $-K_W$ is nef. Furthermore looking at the description of the normal bundle of $E$ in $Y$, $\mathcal{N} _{E/Y}$, in \cite{Cutkosky}, we get 
\[(-K_W)^{3}=(-K_Y)^{3}+a^{3}E^{3}\ge a^{3}E^{3}>0 \]
(note that $(\varphi^{*}K_W+aE)^{3}={K_W}^3+a^{3}E^{3}$ since $\varphi(E)$ is a point). Hence $q(W)=0$ by Kawamata-Viehweg Vanishing Theorem. This derives contradiction.
 
 Hence we shall assume that $\varphi(E)$ is a curve from now on. We know that $C_{0}:=\varphi(E)$ is locally a complete intersection and $\varphi$ is just the blow-up of $W$ along $C_{0}$ (cf. \cite{Cutkosky}). Then $W$ being Cohen-Macaulay implies $E=\mathbb{P}(\mathcal{N}_{{C_0}/W}^{*})$, where $\mathcal{N}_{{C_0}/W}$ denotes the normal bundle of $C_0$ in $W$. Let $\nu:{C_0}'\to C_0$ be the normalization, $V$ the normalized sheaf of $\nu^{*}(\mathcal{N}_{{C_0}/W}^*)$, i.e. $H^{0}(V)\neq0$, but $H^{0}(V\otimes\mathcal{L})=0$ for all line bundles $\mathcal{L}$ with $\deg\mathcal{L}<0$. Furthermore let $\psi:\mathbb{P}(V)(=:E')\to{C_0}'$ be the projection and $F$ its fiber. The tautological line bundle $\mathcal{O}_{\mathbb{P}(V)}(1)$ has a canonical section $C_1$ satisfying ${C_1}^2=c_1(V)(=:-e)$. Put $g$ the genus of ${C_0}'$.
\[
{\xymatrix{
X & Y \supset E \ar[l]_{\!\!\! f} \ar@<-2.8ex>[d]^{\varphi} \ar@<2.5ex>[d] & E' \ar[l]_{\quad\ h} \ar[d]^{\psi} \\
  & W \supset C_{0} & C_{0}' \ar[l]_{\quad\ \nu} }}
\]
We have 
\setcounter{equation}{0}
\begin{equation}
 K_Y=\varphi^{*}K_W+E,         
\end{equation}
and since $\mathcal{N}_{E/Y}=\mathcal{O}_{E}(-1)$
\begin{equation}
 h^{*}(\mathcal{N}_{E/Y})\equiv -C_1+aF  \qquad(a\in\mathbb{Z})
\end{equation}
where $h:E'\to E$ is the base change of $\nu$ by $\varphi$. We can write $V=\nu^{*}(\mathcal{N}_{{C_0}/W}^*)\otimes\mathcal{M}$ for some $\mathcal{M}\in \textup{Pic}(C'_{0})$. Then note that $a=\deg\mathcal{M}$. Moreover
\begin{equation}
 h^{*}(-K_Y|_E)\equiv C_1+bF \qquad(b\in\mathbb{Z})
\end{equation}
from $h^{*}(K_Y|_E)\cdot F=h^{*}(\mathcal{N}_{E/Y})\cdot F=-1.$ Combining (1),(2) and (3), we get 
\begin{align}
 -K_{W}\cdot C_{0} &=-K_{W}\cdot \varphi_{*}h_{*}C_{1}=(-K_{Y}+E)\cdot h_{*}C_{1} \nonumber \\ 
 &=h^{*}(-K_{Y}+E|_E)\cdot C_{1}=a+b,
\end{align}
\begin{equation}
 (-K_{Y})^{2}E=h^{*}(-K_Y|_E)^{2}=-e+2b,
\end{equation}
\begin{equation}
 E^{2}(-K_Y)=h^{*}(\mathcal{N}_{E/Y})\cdot h^{*}(-K_{Y}|_E)=e+a-b,
\end{equation}
\begin{equation}
 (-K_W)^{3}=3{K_Y}^{2}E-3K_{Y}E^{2}+E^{3}=-e+a+3b,
\end{equation}
 and
\[ h^{*}(K_{E})=h^{*}(K_{Y}+E|_E)\equiv -2C_{1}+(a-b)F. \]
On the other hand, since $\omega _{E'}$ is a subsheaf of $h^{*}\omega _E$, the last equation implies 
\[K_{E'}+cF\equiv -2C_{1}+(a-b)F \qquad (c\ge0). \]
Squaring yields
\begin{equation}
 2(1-g)-c=-e-a+b.
\end{equation}

 Note that $-K_{W}\cdot C\ge0$ for any curves $C$ in $W$, $C\ne C_{0}$, and so if $a+b\ge0$, $-K_W$ is nef by (4). Furthermore we know $-e+2b\ge 0$ by (5).

 If $e+a-b>0$, we get 
\[-K_{W}\cdot C_{0}=a+b\ge e+a-b>0 \]
and
\[(-K_{W})^{3}=-e+a+3b\ge e+a-b>0.\]
Hence $-K_W$ is nef and big, so this contradicts $q(W)=q(X)\ne 0$. Consequently we know
\[ e+a-b\le 0.\]
So we obtain from (8) that $2(1-g)\ge c$. Since $c\ge 0$, it is enough to consider the following four cases: \\

 Case 1:\quad $g=0$, $c=0$, Case 2:\quad $g=0$, $c=1$,\\

 Case 3:\quad $g=0$, $c=2$, and Case 4:\quad $g=1$, $c=0$. \\
 
\noindent Before we consider the four cases above, we prove the following claim.\\

 \underline{Claim.} If $g=0$, then $-e+2b>0$.\\
 
\noindent \textit{Proof of the claim}. Assume that $-e+2b=0$. We have   
\[ 0\le h^{*}(-K_{Y}|_E)\cdot C_{1}=-e+b=-\frac{1}{2}e. \] 
Thus $e=0$ by $g=0$ and $K_{Y}\cdot h_{*}C_{1}=0$. Therefore $E'$ is $\mathbb{P}^{1}\times \mathbb{P}^{1}$ and $f\circ h(C_{1}')$ is a point for any curves $C_{1}'\in |C_{1}|$ where $|C_{1}|$ is the complete linear system on $E'$ defined by $C_{1}$. Using (8), we obtain that $a=-2$ in Case 1, $a=-1$ in Case 2 and $a=0$ in Case 3.
\[
{\xymatrix{
X & Y_{0} \ar[l]_{g_{0}} & \cdots \ar[l]_{g_{1}} & Y_{n-1} \ar[l]_{g_{n-1}} & Y_{n}=Y \supset E \ar[l]_{g_{n}} \ar[d]^{\varphi} \ar@<5.2ex>[d] & E'=\mathbb{P}^{1}\times\mathbb{P}^{1} \ar[l]_{\quad\ h} \ar@<-5.2ex>[d]^{\psi} \\
& & & & \ \qquad W \supset C_{0} & C_{0}'=\mathbb{P}^{1}\qquad \ar[l]_{\quad\ \nu} }}
\]
Let $D_i$ be the exceptional divisor of $g_i$ for $i\ge1$. Then there exists $m\ge1$ such that $D_m=g_{m+1}\circ \cdots \circ g_n(E)$ and $\dim g_m(D_m)=1$. We denote the image of $C_1$ on $Y_m$ by $C$ and the strict transform of $D_i$ on $Y$ by $\tilde{D_i}$. We know that $g_m(C)$ is a point and $\dim {g_m}^{-1}(g_m(C))=1$. Thus $C=\mathbb{P}^1$ by $R^1 {g_m}_* \mathcal{O}_{Y_m}=0$. Put ${g_n}^*\cdots {g_{m+1}}^*D_m=E+\sum_{i=m+1}^{n}a_i \tilde{D_i}$ for some $a_i\in \mathbb{Q}_{\ge0}$. Then by (2),
 \[D_m\cdot C=(E+\sum a_i \tilde{D_i})\cdot h(C_1)\ge E\cdot h(C_1)=a.\]  
On the other hand, let $H$ be a general hyperplane section of $Y_{m-1}$. Then the induced morphism $g_{m,H}:H':={g_m}^* H\to H$ is a crepant partial resolution. Now that $C$ is a smooth rational curve,
\[D_m\cdot C=D_m|_{H'}\cdot C=C^2=-2.\]
Therefore we obtain $a=-2$. Moreover we get $c=0$ by (8), so $C_0=\mathbb{P}^1$, $E=\mathbb{P}^{1}\times \mathbb{P}^{1}$ and $\mathcal{N}_{C_0/W}=\mathcal{O}_{\mathbb{P}^1}(-2)\oplus\mathcal{O}_{\mathbb{P}^1}(-2)$ by the definition of $a$ and $c$. Then ${K_Y}^2\equiv 0$ by Proposition 3.11 in \cite{DPS93}. This yields a contradiction because ${K_X}^2 F_{\alpha}>0$ for a general fiber $F_{\alpha}$ of $\alpha$.\\
 
 Now we shall investigate each case.\\
 
\underline{Case 1: $g=0$ and $c=0$}.

 In this case, we know that $C_{0}$ is normal by the definition of $c$, and so $\nu$ and $h$ are isomorphisms. From (8), we obtain
\begin{equation}
 -e-a+b=2.
\end{equation} 
If $-e+2b\ge2$, then
\[(-K_{W})^3=-e+a+3b\ge e+a-b+4=2.\]
Since $-K_{W}$ is neither nef or big, 
\[-K_{W}\cdot C_{0}=a+b<0.\]
The absurdity comes from
\[0=e+a-b+2\le a+b<0.\]
Then the claim above yields
\begin{equation}
-e+2b=1.
\end{equation}
From (9) and (10), $-1=a+b$ and $a=-\frac{1}{2}e-\frac{3}{2}$. Then

\begin{align*}
-1&=-K_{W}\cdot C_{0}=(-K_{Y}+E)\cdot C_{1}\ge E\cdot C_{1}=\mathcal{N}_{E/Y}\cdot C_{1} \nonumber \\
 &=(-C_{1}+aF)\cdot C_{1}=e+a=\frac{1}{2}e-\frac{3}{2}.
\end{align*}
Hence
\[\quad e=0,1.\]
Therefore (10) induces $e=1$. Hence $E$ is a Hirzebruch surface of degree 1 over a smooth rational curve $C_0$. But this implies that $q(Y)=0$ by Proposition 3.5 in \cite{DPS93}.\\

\underline{Case 2: $g=0$ and $c=1$}.
 
 Now $e+a-b=-1$ by (8), and so by the claim above, we get
\[(-K_{W})^3=-e+a+3b\ge e+a-b+2=1.\]
Thus
\[-K_{W}\cdot C_{0}=a+b<0.\]
Then the absurdity comes from
\[-1=e+a-b<a+b<0.\]\\
 
\underline{Case 3: $g=0$ and $c=2$}.

 In this case, $e+a-b=0$ by (8). We can derive contradiction in the similar way as Case 2.\\
 
\underline{Case 4: $g=1$ and $c=0$}.

  In this case, $e+a-b=0$, and $C_0$ and $E$ are smooth by the definition of $c$. If $-e+2b>0$, we obtain that $-K_W$ is nef and big easily. So we know $-e+2b=0$ and \[ 0\le -K_{Y}|_E\cdot C_{1}=-e+b=-\frac{1}{2}e. \] 
Hence $e=0$ or $-1$ because $E$ is an elliptic ruled surface. Combining this with $e=2b$, we know $e=b=0$. Therefore $K_Y\cdot C_{1}=0$. If $\dim f(E)=1$, then there exist infinitely many irreducible curves on $E$ such that each curve is contracted by $f$. Take such an irreducible curve $C$ on $E$. Then by (3), $C$ is numerically proportional to $C_1$ in $N^1(E)$. If $\dim f^{-1}(f(C))=1$, then $C$ is a smooth rational curve. But this derives contradiction because $\varphi(C)$ is an elliptic curve by $g=1$. Thus there exists a prime divisor $D$ such that $f(D)=f(C)$. This is also absurd because $f$ never contracts infinitely many divisors. Therefore we get that $\dim f(E)=2$. Let $H$ be a general hyperplane section of $f(E)$. Then $H':=(f|_{E})^*H$ is an irreducible curve on $E$ which dose not intersect $C_1$. So we know that $H'\equiv dC_1$ for some $d>0$ (cf. Chap.V, 2.20 in \cite{Ha77}). But since $f(H')$ is not a point, this lies in contradiction.
\end{proof}

\section{Proof of $\boldsymbol{\alb(X)\ne2}$}
\indent

The proof of the following result is similar to the one of Theorem 1 in \cite{Zhang}, which treats the smooth case in arbitrary dimension.

\begin{prop}\label{Zhang}
Let $X$ be a canonical 3-fold with nef anticanonical divisor and
$\alpha :X\to \Alb(X)$ its Albanese morphism. Assume
$\alb(X)=2$. Then $\alpha$ is an algebraic fiber space.  
\end{prop}
\begin{proof}
It is enough to consider the case that $X$ is a terminal 3-fold by terminalization. Since a general hyperplane section of a terminal 3-fold is smooth, we can use the proof of \cite{Zhang}.
\end{proof}

\begin{prop} \textbf{\textup{\cite{Kollar}}} \label{Kollar}
Let $\varphi :X\to Y$ be a surjective morphism between smooth
varieties. Assume that the general fiber $F$ is connected and
$h^{i}(F,\mathcal{O}_F)=0$ for $i>0$. Then
$h^{i}(X,\mathcal{O}_X)=h^{i}(Y,\mathcal{O}_Y)$ for all $i$.
\end{prop}

\begin{prop} \textbf{\textup{\cite{PS}}}\label{lem:PS2}
Let $X$ be a smooth projective 3-fold with nef $-K_{X}$ and $q(X)=2$. We stand $\alpha$ for the Albanese morphism of $X$ and $\varphi :Y\to X$ for a divisorial contraction from a terminal 3-fold $Y$. Assume that $\alpha$ is a smooth morphism. Then $-K_{Y}$ is not almost nef. 
\end{prop}

\begin{prop}\label{alb=2}
Let $X$ be a canonical 3-fold with strictly nef $-K_X$. Then $\alb(X)\ne2$.
\end{prop}
 
\begin{proof}
Assume that $\alb(X)=2$. Write $A:=\Alb(X)$. Now that the Albanese morphism $\alpha$ is an algebraic fiber space (Proposition \ref{Zhang}), a general fiber of $\alpha$ is $\mathbb{P}^{1}$. Using Proposition \ref{Kollar}, we know $h^{i}(\mathcal{O}_X)=h^{i}(\mathcal{O}_A)$.

 Let $f:Y\to X$ be a $\mathbb{Q}$-factorial terminalization as in the proof of Proposition \ref{alb=1}. Since $K_Y$ is not nef, there exists an extremal contraction $\varphi :Y\to W$ and we know $\dim W\ge2$ because Proposition \ref{Zhang} says in paticular that $\alpha$ is surjective.\\ 

\textbf{\underline{Case: $\dim W=2$}}
 
Now $h^{0}(K_W)=h^{2}(\mathcal{O}_W)=h^{2}(\mathcal{O}_X)=1$, and so $K_{W}\cdot H\ge0$ for some ample divisor $H$ on $Y$. On the other hand, since $-(4K_W+\Delta)$ is almost nef (Proposition \ref{prop:PS}), $-(4K_W+\Delta)\cdot H\ge0$. Thus $K_W\equiv0$. So we know that $W$ is an abelian surface and hence $\beta$ is an isomorphism where $\beta:W\to \Alb(W)(=A)$ is the Albanese morphism.

 If $X$ is not $\mathbb{Q}$-factorial terminal, there exists a curve $C$ on $Y$ such that $f(C)$ is a point. On the other hand, $K_{Y}\cdot C=0$ implies $\varphi(C)$ is a curve. This contradicts $\varphi=\alpha\circ f$. Hence $X$ is a $\mathbb{Q}$-factorial terminal 3-fold and we can apply Proposition \ref{prop:PS1}, so we know that $X$ is smooth. Then $-K_X$ is ample by Theorem \ref{Serrano3} (ii). The ampleness of $-K_X$ yields a cotradiction with $\alb(X)=2$. \\   

\textbf{\underline{Case: $\dim W=3$}} 
 
The following argument is in pages 493-494 of \cite{PS}. By Mori Theory, there exists a sequence  
\[ Y\to W=:Y_1\to\dots\to Y_{n-1} \to Y_n \to Z. \] 
where each morphism $\varphi_i:Y_i\to Y_{i+1}$ is a divisorial contraction or a flip for $i\ge1$ and $\psi:Y_n\to Z$ is a Mori fiber space to a surface $Z$. As in the case $\dim W=2$ above, we know that $Z$ is an abelian surface, so we get readily that $\psi$ is the Albanese morphism of $Y_n$. So by Proposition \ref{prop:PS1}, $\psi$ is a smooth morphism and $-K_{Y_n}$ is nef. Thus $\varphi_{n-1}$ is divisorial. But this is absurd with  Proposition \ref{prop:PS} and Proposition \ref{lem:PS2}.
\end{proof}
 
\section{Proof of the Main Theorem}

\indent

We first prepare several lemmas in order to prove Theorem \ref{q=0}.
\begin{lem}\label{H2<2}
Let $X$ be a canonical 3-fold with $-K_{X}$ nef and $\kappa(X)=-\infty$. Then $h^{2}(\mathcal{O}_X)<2$.
\end{lem}

\begin{proof}
Taking a $\mathbb{Q}$-factorial terminalization, we may assume that $X$ has at most $\mathbb{Q}$-factorial terminal singularities. By Mori theory, $X$ is birational to a $\mathbb{Q}$-factorial terminal 3-fold $W$ such that $W$ has a Mori fiber space structure $\varphi :W\to Z$, that is, $\varphi$ is an extremal contraction with $\dim W>\dim Z$. If $h^{2}(\mathcal{O}_X)\geq 2$, then $Z$ is a projective surface and  $h^{0}(K_{Z})=h^{2}(\mathcal{O}_Z)=h^{2}(\mathcal{O}_X)\geq 2$. So $K_{Z}\cdot H\geq 1$ for some ample divisor $H$. On the other hand, since $-(4K_{Z}+\Delta)\cdot H\geq 0$ (Proposition \ref{prop:PS}), $-K_{Z} \cdot H\geq 0$. This is absurd.
\end{proof}

\begin{lem} \label{finite}
Let $X$ be a canonical 3-fold such that $-K_{X}$ is strictly nef. Assume $q(X)=0$. Then $\pi^{alg}_1(X)$ is finite.
\end{lem}

\begin{proof}
Assume that $\pi^{alg}_1(X)$ is infinite. Then there exists an infinite tower of normal finite \'etale Galois covers
\[ \dots\to X_2 \to X_1 \to X_0:=X. \]
We know $\chi(\mathcal{O}_{X_m})\gg 0$ by $\chi(\mathcal{O}_X)>0$, hence $h^{2}(\mathcal{O}_{X_m})\gg 0$ for a sufficiently large $m$. Then $-K_{X_m}$ being strictly nef contradicts Lemma \ref{H2<2}.
\end{proof}

\begin{lem}\textbf{\textup{\cite{Miyaoka}}}\label{lem:Miya}
Let $S$ be a smooth projective surface with finite algebraic fundamental group $\pi^{alg}_1 (S)$. Suppose we have a nef divisor $D$ on $S$ such that $D^{2}=0$. Then either $H^{1}(S,-D)=0$, or there exists a positive integer $n$ such that $H^{0}(S,nD)\ne 0$.
\end{lem} 
The following lemma is the generalization of Lemma 3.8 in \cite{Serrano}. 
\begin{lem}\label{K2=0}
Let $X$ be a canonical 3-fold with strictly nef $-K_X$. Let $U$ stand for the locus of cDV points of $X$. Assume that ${K_X}^2 \equiv 0$ and $\pi^{alg}_1(U)$ is finite. Then $\chi(\mathcal{O}_X) \le 0$.  
\end{lem}

\begin{proof}
Because $\pi^{alg}_1(U)=\pi^{alg}_1(S)$ for a general hyperplane section $S$ by Th\'{e}or\`{e}me 3.10, page 123 of \cite{SGA2}, $\pi^{alg}_1(S)$ is finite. Assume that $\chi(\mathcal{O}_X)> 0$ below. Write $r:=\mathrm{index}(X)$. Let $f:Y\to X$ be a terminalization and set $T:=f^*S$. Since $-K_{Y}$ is nef and not big, $c_1(Y)c_2(Y)\ge 0$ by Proposition \ref{KeMM} and Remark \ref{KeMM rem}. If  $h^{0}(-rK_{Y})=h^{0}(-rK_{X})\ne 0$, then $-K_{X}$ is ample by Theorem \ref{Serrano3} (iii) and so it contradicts ${K_X}^{2}\equiv 0$. Hence Riemann-Roch Theorem (cf. \cite{Reid87}) says that
 \[ \chi(-rK_{Y})=\frac{1}{12} (-rK_Y)c_2(Y)+\chi(\mathcal{O}_Y)>0,\] 
 which implies $h^2 (-rK_Y)>0$. Let us consider the exact sequence
\begin{align*}
\cdots&\to H^{1}(-rK_Y+T) \to H^{1}(-rK_{Y}+T|_{T}) \to H^{2}(-rK_{Y}) \to \\
       &\to H^{2}(-rK_{Y}+T)\to\cdots .
\end{align*}
Exchange $S$ for its multiple if necessary. Then we may assume $h^{i}(-rK_{Y}+T)=h^{i}(-rK_{X}+S)=0$ for $i>0$. Therefore
\[h^{1}((r+1)K_Y|_T)=h^{1}(-(r+1)K_Y|_T+K_T)=h^{1}(-rK_Y+T|_T)\ne 0.\] 
Now that $\pi^{alg}_1(T)=\pi^{alg}_1(S)$ (cf. Lemma 3.1 in \cite{SBW}) and $T$ is smooth, we can apply Lemma \ref{lem:Miya} and we get $h^{0}(-nK_X|_S)=h^{0}(-nK_Y|_T)\ne 0$ for some $n>0$. Hence
\[ 0<-K_X\cdot (-nK_X|_S)=n(-K_X)^2\cdot S=0. \]
This is a contradiction. So $\chi(\mathcal{O}_X) \le 0$.  
\end{proof}

\begin{rem} \label{K3>0}
Let $X$ be a canonical 3-fold with strictly nef $-K_X$. Assume $(-K_{X})^3>0$. Then we know that $-K_X$ is ample. In fact, the complete linear system $|-nK_{X}|$ is free for $n\gg 0$ by the Base Point Free Theorem (cf. \cite{KMM}). Thus it yields a morphism $\Phi: X\to \mathbb{P}^N$. If there exists a curve $C$ on $X$ such that $\Phi (C)$ is a point, then $-nK_{X}\cdot C=\Phi ^*H\cdot C=0$ for a hyperplane $H$ in $\mathbb{P}^N$, which contradicts the strictly nefness of $-K_X$. Therefore $\Phi$ is finite, and we obtain that $-K_X$ is ample.      
\end{rem} 

\begin{lem}\label{K2>0}
Let $X$ be a canonical 3-fold with strictly nef $-K_X$. Assume that $q(X)=0$ and ${K_X}^2 \not\equiv 0$. Then $-K_X$ is ample.
\end{lem}

\begin{proof}
Let $S$ be a sufficiently ample general hyperplane section and $r:= \mathrm{index}(X)$. Consider the exact sequence
\[ \begin{CD}
 @>>> H^{1}(-rK_{X}+S|_{S}) @>>> H^{2}(-rK_{X}) @>>> H^{2}(-rK_{X}+S) @>>> . \\
@. @|  @.  @|  @. \\
@. H^{1}(-(r+1)K_{X}|_{S}+K_{S}) @.  @.  0  @. \\
@. @|  @.  @.  @. \\
@. 0  @. @. @.
\end{CD} \]
Therefore $h^{2}(-rK_X)=0$. Let $f:Y\to X$ be a terminalization. If $(-K_X)^3=(-K_{Y})^3=0$, then $c_1(Y)c_2(Y) \ge 0$ by Proposition \ref{KeMM} and Remark \ref{KeMM rem}. Furthermore by $\chi(\mathcal{O}_Y)=\chi(\mathcal{O}_X)>0$,
 \[ \chi(-rK_{Y})=\frac{1}{12} (-rK_Y)c_2(Y)+\chi(\mathcal{O}_Y)>0.\] 
 So $h^{0}(-rK_X) \ne 0$. This implies $-K_X$ is ample by Theorem \ref{Serrano3} (iii). On the other hand, $-K_X$ is ample when $(-K_{X})^3 >0$ (Remark \ref{K3>0}). This completes the proof.
\end{proof}

\begin{rem} \label{q==0}
 Let $X$ be a canonical 3-fold with strictly nef $-K_X$. Then we have $q(X)=0$. In fact, because $\alpha (X)$ is a subvariety of $\Alb (X)$, $\kappa(\alpha(X))\ge 0$. So $\alb(X)=3$ induces a contradiction with $\kappa(X)=-\infty$. Thus Proposition \ref{alb=1} and Proposition \ref{alb=2} deduce that $\alb(X)=0$, which is equivalent to $q(X)=0$. 
\end{rem} 

We are now in the position to prove the main theorem of this paper. The proof is based on the Miyaoka's idea in \cite{Miyaoka}. 
\begin{thm} \label{q=0}
Let $X$ be a canonical 3-fold with strictly nef $-K_X$. Then $-K_{X}$ is ample. 
\end{thm}
\begin{proof}
 Note that $q(X)=0$ by Remark \ref{q==0}. Let $U$ be the locus of cDV points of $X$. If $\pi^{alg}_1(U)$ is finite, combining Lemma \ref{K2=0} and Lemma \ref{K2>0}, we know that $-K_X$ is ample because $\chi(\mathcal{O}_X)>0$.

 If $\pi^{alg}_1(U)$ is infinite, then there exists an infinite tower of normal finite Galois covers
\[ \dots\to X_2 \to X_1 \to X_0:=X \]
such that $f_{n}|_{U_n} :U_n\to U$ is \'etale where $f_n$ is the morphism from $X_n$ to $X$ and $U_n:=f^{-1}_{n} (U)$. Let $x$ be a non-cDV point in $X$ and $g_n$ the morphism from $X_n$ to $X_{n-1}$. Since $\pi^{alg}_1(V^o)$ is a finite group (Theorem 3.6 in \cite{SBW}) where $V$ is a sufficiently small analytic neighbourhood of $x$ and $V^o$ stands for the smooth locus of $V$ (we use $^o$ below as the same way), there exists a positive integer $m_x$ such that $f_{n+1}^{-1}(V)^o$ is disjoint union of copies of $f_{n}^{-1}(V)^o$ for all $n\ge m_x$. Now that $X_{n}$ is normal, for all $y\in f_{n}^{-1}(x)$, $f_{n+1}^{-1}(V)^o\cup g_{n+1}^{-1}(y)$ is also disjoint union of copies of $f_{n}^{-1}(V)^o\cup y$ (see Th\'{e}or\`{e}me 3.6, page 38 of \cite{SGA2}) and so $g_n$ is \'etale over $y$. $X$ has only finitely many non-cDV points, so $g_{n}$ is \'etale and $\pi^{alg}_1(X_n)$ is infinite for all sufficiently large $n$. Now that $X_n$ is a canonical 3-fold with strictly nef $-K_{X!
_n}$, we derive a contradiction from Lemma \ref{finite}.
\end{proof} 

From Lemma \ref{H2<2}, we also obtain the following result: 
\begin{thm}\textbf{\textup{(cf. \cite{Zhang})}}
Let X be a canonical 3-fold with nef $-K_X$. Then $q(X)\le3$. Furthermore $q(X)=3$ if and only if $X$ is an abelian variety.
\end{thm}

\begin{proof}
Assume that $\kappa(X)\ge0$. Then we know readily that $\kappa(X)=0$. Let $f:Y\to X$ be a resolution of singularities. Because the Albanese morphism of $Y$ is an algebraic fiber space (\cite{K81}), so is the one of $X$, $\alpha: X\to\Alb (X)=:A$ and $q(X)\le3$. Consider the case that the equality holds. Then we can regard terminalization of $X$ as a resolution of singularities $f$ (Corollary 8.4 in \cite{K85}). So $Y$ is an abelian variety by \cite{Zhang}. Hence $f$ and $\alpha$ are isomorphic.

 If $\kappa(X)=-\infty$, it is enough to consider the case that $-K_X$ is not big. According to Proposition \ref{KeMM} and Remark \ref{KeMM rem}, $\chi(\mathcal{O}_X)\ge0$ so we get $q(X)<3$ by Lemma \ref{H2<2}.  
\end{proof}

\textsc{Department of Mathematical Sciences, University of Tokyo, Komaba, Meguro, Tokyo 153-8914, Japan} \\
\texttt{hokuto@ms318sun.ms.u-tokyo.ac.jp}


\begin{thebibliography}{85}
\bibitem{Benveniste}
X. Benveniste, Sur le cone des 1-cycles effectifs en dimension 3,
Math. Ann. 272 (1985), 257-265.
\bibitem{Cutkosky}
S. Cutkosky, Elementary contractions of Gorenstein threefolds,
Math. Ann. 280 (1988), 521-525. 
\bibitem{DPS93}
J. P. Demailly, Th. Peternell, M. Schneider, K\"ahler manifolds with numerically effective Ricci class, Comp. Math. 89 (1993), 217-240. 
\bibitem{Elkik}
R. Elkik, Rationalit\'{e} des singularit\'{e}s canoniques, Invent. Math. 64 (1981), 1-6. 
\bibitem{Fu97}
W. Fulton, Intersection Theory, 2nd. ed. A Series of Modern Surveys in Matematics, vol. 2, Springer-Verlag, (1997).
\bibitem{SGA2}
A. Grothendieck, Cohomologie locale des faisceaux coh\'{e}rents et th\'{e}or\`{e}me de Lefschetz locaux et globaux (SGA2), Masson, North-Holland, Paris-Amsterdam, (1962).
\bibitem{Ha70}
R. Hartshorne, Ample subvarieties of algebraic varieties, Lect. Notes Math. 156, Springer-Verlag, (1970).
\bibitem{Ha77}
R. Hartshorne, Algebraic geometry, Grad. Texts. in Math. 52, Springer-Verlag, (1977).
\bibitem{K81}
Y. Kawamata, Characterization of Abelian varieties, Comp. Math 43 (1981), 253-276.
\bibitem{K85}
Y. Kawamata, Minimal models and the Kodaira dimension of algebraic fiber spaces, J. reine angew. Math. 363 (1985), 1-46.
\bibitem{K86}
Y. Kawamata, On the plurigenera of minimal algebraic 3-folds with $K\equiv 0$, Math. Ann. 275 (1986), 539-546.
\bibitem{K88}
Y. Kawamata, Crepant blowing-up of 3-dimensional canonical singularities and its application to degenerations of surfaces, Ann. of Math. 127 (1988), 93-163.
\bibitem{KMM}
Y. Kawamata, K. Matsuda, K. Matsuki, Introduction to the minimal model 
problem, Adv. Stud. Pure Math. 10, Alg. Geom. Sendai 1985 (T.Oda,
ed.), Kinokuniya, Tokyo (1987), 283-360.
\bibitem{KeMM}
S. Keel, K. Matsuki, J. McKernan, Log abundance theorem for threefolds, Duke 
Math. J. 75 (1994), 99-119.
\bibitem{Kleiman}
S. Kleiman, Toward a numerical theory of ampleness, Ann. of Math. 84 (1966), 293-344.   
\bibitem{Kollar}
J. Koll\'{a}r, Higher direct images of dualizing sheaves I,
Ann. of Math. 123 (1986), 11-42. 
\bibitem{Kollar89}
J. Koll\'{a}r, Flops, Nagoya Math. J. 113 (1989), 15-36.
\bibitem{KM}
J. Koll\'{a}r, S. Mori, Birational geometry of algebraic varieties, Cambridge Tracts in Math, vol. 134, (1998).
\bibitem{Maeda}
H. Maeda, A criterion for a smooth surface to be Del Pezzo,
Math. Proc. Cambridge Phil. Soc. 113 (1993), 1-3.
\bibitem{Matsuki}
K. Matsuki, A criterion for the canonical bundle of a 3-fold to be
ample, Math. Ann. 276 (1987), 557-564. 
\bibitem{Miyaoka}
Y. Miyaoka, On the Kodaira dimension of minimal threefolds, Math. Ann. 281 (1988), 325-332.
\bibitem{PS}
Th. Peternell, F. Serrano, Threefolds with nef anticanonical bundles, Collect. Math. 49 (1998), 465-517.
\bibitem{Reid83}
M. Reid, Minimal models of canonical 3-folds, Adv. Stud. Pure Math, vol.1, Kinokuniya, Tokyo (1983), 131-180. 
\bibitem{Reid87} 
M. Reid, Young person's guide to canonical
singurarities, Proc. Symp. Pure Math. 46 (1987), 345-414.
\bibitem{Serrano}
F. Serrano, Strictly nef divisors and Fano threefolds, J. reine
angew. Math. 464 (1995), 187-206. 
\bibitem{SBW}
N. I. Shepherd-Barron, P. M. H. Wilson, Singular threefolds with numerically trivial first and second chern classes, J. Alg. Geom. 3 (1994), 265-281.
\bibitem{Zhang}
Q. Zhang, On projective manifolds with nef anticanonical bundles,
J. reine angew. Math. 478 (1996), 57-60.
\end{thebibliography}
\end{document}